\documentclass{aims} 
  \usepackage{paralist}
  \usepackage{graphics} 
  \usepackage{epsfig} 
\usepackage{graphicx}  
\usepackage[colorlinks=true]{hyperref}
\hypersetup{urlcolor=blue, citecolor=red}
\newcommand{\sg}{\mbox{sgn}}
\usepackage{bm}
\allowdisplaybreaks

\def\currentvolume{x}
 \def\currentissue{x}
  \def\currentyear{xx}
   \def\currentmonth{xx}
    \def\ppages{x--x}
      \def\DOI{x}

\newtheorem{theorem}{Theorem}[section]

\theoremstyle{definition}
\newtheorem{definition}[theorem]{Definition}

\title[Globally stable periodic orbits]
{A variety of globally stable periodic orbits in permutation binary neural networks}

\author[Mikito Onuki, Kento Saka and Toshimichi Saito]{}


 \email{tsaito@hosei.ac.jp}

\thanks{$^*$  Corresponding author: Toshimichi Saito}

\begin{document}
\maketitle

\centerline{\scshape Mikito Onuki,  Kento Saka and Toshimichi Saito$^*$}
\medskip
{\footnotesize
 \centerline{Department of Electrical and Electronic Engineering}
   \centerline{ HOSEI University, Japan}
} 

%

\bigskip



\begin{abstract}
The permutation binary neural networks are characterized by global permutation connections and local binary connections. 
Although the parameter space is not large, the networks exhibit various binary periodic orbits.  
Since analysis of all the periodic orbits is not easy, we focus on 
globally stable binary periodic orbits such that almost all initial points fall into the orbits. 
For efficient analysis, we define the standard permutation connection that represents multiple equivalent permutation connections. 
Applying the brute force attack to 7-dimensional networks, we present the main result: 
a list of standard permutation connections for all the globally stable periodic orbits.  
These results will be developed into detailed analysis of the networks and its engineering applications.  
\end{abstract}

\section{Introduction}
\label{intro}
Discrete-time recurrent neural networks (DT-RNNs) 
are analog dynamical systems characterized by 
real valued connection parameters and nonlinear activation functions (e.g., sigmoid function) \cite{rnn1} \cite{rnn2} \cite{rnn3} \cite{rnn4}.  
The dynamics is described by autonomous difference equations of real state variables. 
Depending on the parameters, the DT-RNNs exhibit various periodic orbits, chaos \cite{ott}, and related bifurcation phenomena. 
The real/potential applications include
associative memories \cite{rnn1}, 
combinatorial optimization problems solvers \cite{tsp}, and 
time-series approximation/prediction in reservoir computing \cite{rc1} \cite{rc2} \cite{rc3}. 
The DT-RNNs are important systems in both 
basic study of nonlinear dynamics and engineering applications. 
However, analysis of the dynamics is hard because of 
huge parameter space and complexity of the nonlinear phenomena. 
Stability analysis of various periodic orbits is not easy. 

The three-layer dynamic binary neural networks (DBNNs \cite{koyama} \cite{anzai}) are digital dynamical systems 
characterized by ternary valued connection parameters and the signum activation function. 
The dynamics is described by autonomous difference equations of binary state variables. 
Since the state space consists of a finite number of binary variables, the DBNNs cannot generate chaos \cite{ott}. 
However, depending on the parameters and initial conditions, the DBNNs can generate various periodic orbits of binary vectors 
(binary periodic orbits, ab. BPOs). 
As compared with the DT-RNNs, the DBNNs bring benefits to hardware implementation. 
An FPGA based hardware prototype and its application to hexapod walking robots can be found in \cite{takumi}.  
We have presented a parameter setting method that guarantees storage and stability of desired BPOs \cite{koyama}. 
However, as period of a BPO increases, the number of hidden neurons increases:  
parameter space becomes wider and analysis becomes harder. 
In the hardware, power consumption becomes larger. 
In order to realize efficient analysis and synthesis, reduction of the parameter space is inevitable. 

Simplifying connection parameters of the DBNNs, the permutation binary neural networks (PBNNs \cite{taka}) are constructed. 
The PBNNs are characterized by two kinds of connections. 
The first one is local binary connection between input and hidden layers. 
It is defined by a signum-type neuron from three binary inputs to one binary output. 
The second one is global one-to-one connection between hidden and output layers. 
It is defined by a permutation operator. 
The parameter space of the PBNNs is much smaller than that of DBNNs. 
Depending on the permutation connections, the PBNNs generate various BPOs. 
Co-existence of BPOs is possible and 
the PBNN exhibits one of the BPOs depending on initial condition. 

Since analysis of multiple BPOs are not easy, 
this paper focuses on globally stable binary periodic orbits (GBPOs) 
such that almost all initial points fall into the GBPOs. 
As a fundamental concept, we define the standard permutation connection  
that represents multiple equivalent permutation connections. 
Applying the brute force attack to all the 7-dimensional PBNNs, 
we present the main result:  
a list of the standard permutation connections for all the GBPOs. 
These results provide basic information to realize more detailed analysis of PBNNs and its applications. 
Real/potential engineering applications of the GBPOs include 
time-series approximation/prediction \cite{rc2} \cite{nara} \cite{uchida}, 
control signals of switching power converters \cite{pe1} \cite{pe2} \cite{pe3}, 
control signals of walking robots \cite{takumi} \cite{cpg}, and
error correcting codes \cite{error}. 
The approximate/control signals can be globally stable and robust. 
As novelty of this paper, it should be noted that 
this is the first paper of the GBPOs and standard permutation connections.

\section{Permutation binary neural networks and binary periodic orbits}
This section introduces 
the 3-layer dynamic binary neural networks (DBNNs, \cite{koyama}) and 
the permutation binary neural networks (PBNNs, \cite{taka}). 
After overview of BPOs, we show the objective problem.   

\subsection{Dynamics binary neural network}
\begin{figure}[tb]
\centering
\includegraphics[width=0.35\columnwidth]{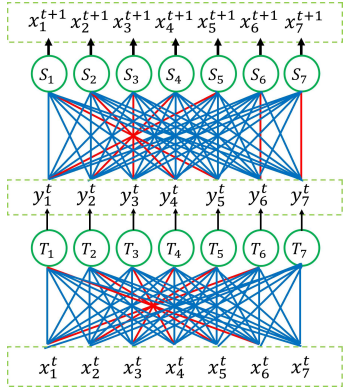}
\caption{Dynamic binary neural network (DBNN)   
Red and blue branches denote positive and negative connections, respectively. 
}
\label{fg1}
\end{figure}

The DBNNs are recurrent-type 3-layer networks characterized by 
ternary connection parameters and signum activation function.  
The dynamics is described by the following autonomous difference equation of $N$-dimensional binary state variables:   
\begin{equation}
\begin{array}{l}
\displaystyle x_i^{t+1} = \sg \left(\sum_{j=1}^M c_{ij} y_j^t  +  S_i \right), 
\ y_j^t = \sg \left(\sum_{i=1}^N w_{ji} x_i^t - T_j \right)\\
\sg(x) = \left\{
\begin{array}{lll}
+1 & \mbox{if } x \ge 0, & \ i \in \{1, \cdots, N \} \\
-1 & \mbox{if } x < 0,   & \ j \in \{1, \cdots, M \}
\end{array}\right. 
\end{array}
\label{3dbnn}
\end{equation}
where 
$x_i^t \in \{-1, +1\} \equiv \bm{B}$ is the $i$-th binary state variable at discrete time $t$ and 
$y_j^t \in \bm{B}$ is the $j$-th binary hidden variable. 
As shown in Fig. \ref{fg1},  the binary variables $x_i^t$, $y_j^t$, and $x_i^{t+1}$ are located in input, hidden, and output layers, respectively. 
The $M$ hidden neurons transform $x_i^t$ into $y_j^t$ through hidden ternary connections ($w_{ji} \in \{-1, 0, +1\}$). 
The $N$ output neurons transform $y_j^t$ into $x_i^{t+1}$ through output ternary connections ($c_{ij} \in \{-1,  0, +1\}$). 
The threshold parameters $S_i$ and $T_j$ are integers. 
The output $x_i^{t+1}$ is fed back to the input layer and the DBNNs generate various BPOs. 
Ref. \cite{koyama} gives a theoretical result of parameter condition 
that guarantees storage and stability of desired BPOs. 
However, as period of a BPO increases, the number of hidden neurons increases. 
For example, $p$ hidden neurons are required for storage of a BPO with period $p$.  
As $p$ increases, the parameter space becomes larger and analysis/implementation becomes harder.

\subsection{Permutation binary neural networks}

\begin{figure}[b!]
  \centering
  \includegraphics[width=0.8\columnwidth]{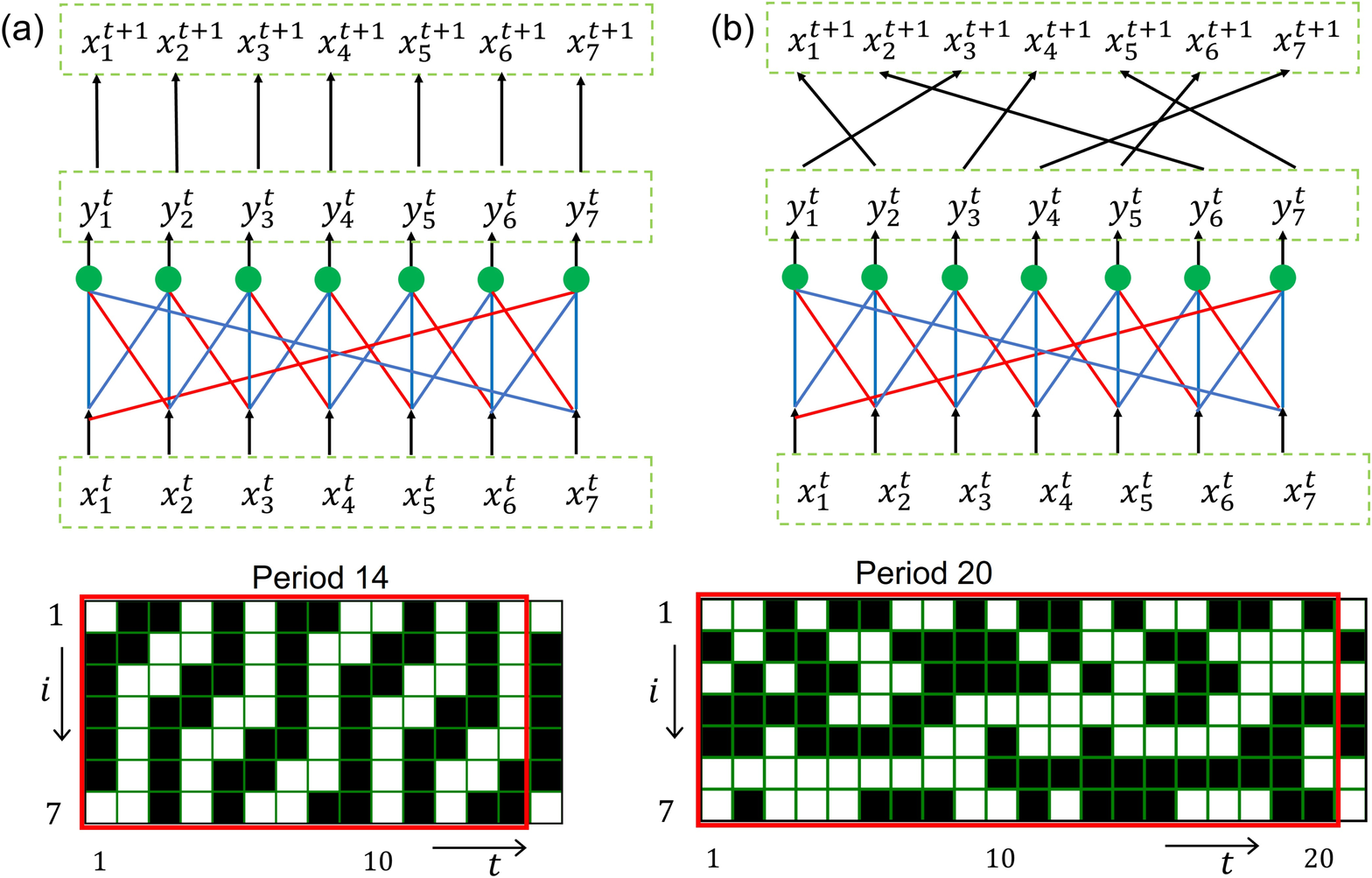}
  \caption{Examples of PBNNs and BPOs for CN1, $N=7$. 
Red and blue branches denote positive and negative local binary connections, respectively. 
black branches correspond to global permutation connections. 
White and black squares in spatiotemporal patterns denote $x_i^t = +1$ and $x_i^t=-1$, respectively. 
(a) $P(1234567)$. (b) $P(2613754)$.
  }
\label{fg2}

\vspace*{5mm}

\centering
\includegraphics[width=0.7\columnwidth]{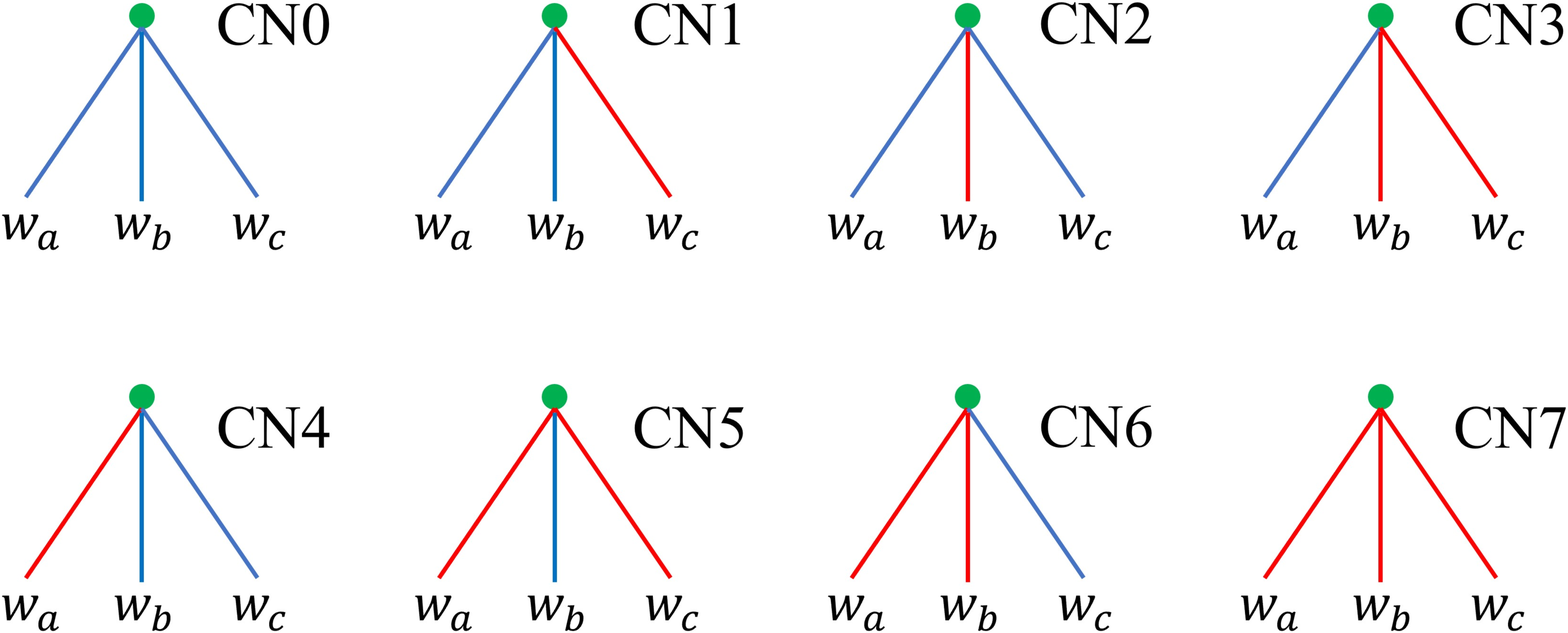}
\caption{8 local binary connections. 
}
\label{fg3}
\end{figure}

The PBNNs are described by the following autonomous difference equation of $N$-dimensional binary state variables:   
\begin{equation}
\begin{array}{l}
x_i^{t+1} = y_{\sigma(i)}^t, \ 
y_i^t = \sg \left(w_a x_{i-1}^t + w_b x_i^t + w_c x_{i+1}^t \right)\\
\sigma = \left(
    \begin{array}{cccc}
      1         & 2         & \cdots & N \\
      \sigma(1) & \sigma(2) & \cdots & \sigma(N)
    \end{array}
  \right) \ i \in \{1, \cdots, N\}, N \ge 3
\end{array}
\label{pbnn}
\end{equation}
where 
$x_0^t \equiv x_N^t$ and $x_{N+1}^t \equiv x_1^t$ for ring-type connection as shown in Fig. \ref{fg2}. 
As a binary state vector $\bm{x}^t \equiv (x_1^t, \cdots, x_N^t) \in \bm{B}^N$ is input at time $t$, 
the $\bm{x}^t$ is transformed into the binary hidden state vector $\bm{y}^t \equiv (y_1^t, \cdots, y_N^t) \in \bm{B}^N$ 
through hidden neurons with local binary connections. 
All the hidden neurons have the same characteristics: 
the signum activation function from three binary inputs to one binary output with local binary connection parameters $(w_a, w_b, w_c) \in \bm{B}^3$. 
The $\bm{y}^t$ is transformed into $\bm{x}^{t+1}$ through one-to-one global permutation connection defined by the permutation $\sigma$. 
The output vector $\bm{x}^{t+1}$ is fed back to the input and the PBNNs generate sequences of binary vectors. 
In comparison with the DBNNs, 
the hidden connections $w_{ij}$ are replaced with the local binary connections and 
the output connections $c_{ij}$ are replaced with the global permutation connections. 
As shown in Fig. \ref{fg3}, the local binary connections are identified by connection numbers:
\[
\begin{array}{ccc}
\mbox{CN0}: \bm{w}_l=(-1, -1, -1) & 
\mbox{CN1}: \bm{w}_l=(-1, -1, +1) & 
\mbox{CN2}: \bm{w}_l=(-1, +1, -1) \\
\mbox{CN3}: \bm{w}_l=(-1, +1, +1) & 
\mbox{CN4}: \bm{w}_l=(+1, -1, -1) &
\mbox{CN5}: \bm{w}_l=(+1, -1, +1) \\ 
\mbox{CN6}: \bm{w}_l=(+1, +1, -1) &
\mbox{CN7}: \bm{w}_l=(+1, +1, +1) 
\end{array}
\]
where $\bm{w}_l \equiv (w_a, w_b, w_c)$. 
Since CN1 (respectively, CN3) coincides with CN4 (respectively, CN6) by replacement 
$x_i \rightarrow x_{N-i+1}$ for $i \in \{1, \cdots, N \}$, 
we consider 6 connection numbers without CN4 and CN6 hereafter. 
The global permutation connections are identified by 
\[
\mbox{Permutation ID: }  P(\sigma(1) \cdots \sigma(N)). 
\]
Fig. \ref{fg2} shows examples of 7-dimensional PBNNs for CN1. 
For identity permutation P(123456), the PBNN exhibits a BPO with period 14. 
Applying permutation P(2613754), the PBNN exhibits a BPO with longer period 20. 
In the DBNN, 20 hidden neurons are necessary for period 20. 

\subsection{Objective problem}

In order to visualize the dynamics, we have introduced the digital return map (Dmap). 
The domain $\bm{B}^N$ of the PBNNs is equivalent to a set of $2^N$ points $L_N \equiv \{ C_1, \cdots, C_{2^N} \}$,
i.e., $C_1 \equiv (-1, \cdots, -1)$, $C_2 \equiv (+1. -1. \cdots, -1)$, $\cdots$, $C_{2^N} \equiv (+1, \cdots, +1)$.
The dynamics of a PBNN can be integrated into
\begin{equation}
\mbox{Dmap: }   \bm{x}^{t+1} = f(\bm{x}^t), \  \bm{x}^t \in \bm{B}^N \equiv L_D
\end{equation}
where an $N$-dimensional binary vector $\bm{x}^t$ is denoted by a point $C_i$ in the Dmap. 

\begin{figure}[tb]
\centering
\includegraphics[width=1.0\columnwidth]{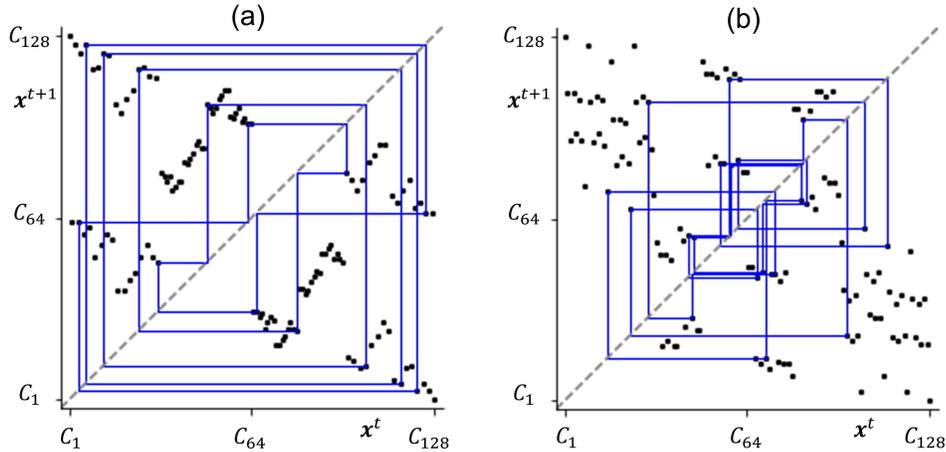}
\caption{Dmap examples (black points) and BPOs (blue orbits) for CN1, $N=7$. 
(a) $P(1234567)$ (the PBNN is Fig. \ref{fg2} (a)), BPO with period 14. 
(b) $P(2613754)$ (the PBNN is Fig. \ref{fg2} (b)), BPO with period 20. 
}
\label{fg4}
\end{figure}

\begin{definition}
\label{def1}
A point $\bm{z}_p \in L_D$ is said to be a binary periodic point (BPP) with period $p$
if $f^p(\bm{z}_p) = \bm{z}_p$ and $f(\bm{z}_p)$ to $f^p(\bm{z}_p)$ are all different
where $f^k$ is the $k$-fold composition of $f$. 
A sequence of the BPPs, $\{ f(\bm{z}_p), \cdots, f^p(\bm{z}_p) \}$,
is said to be a BPO with period $p$.
A point $\bm{z}_e$ is said to be an eventually periodic point (EPP)
if $\bm{z}_e$ is not a BPP but falls into a BPO, i.e.,
there exists some positive integer $l$ such that $f^l(\bm{z}_e)$ is a BPP. 
The BPO in the Dmap is equivalent to the BPO in spatiotemporal pattern from the PBNN.  
\end{definition}
Fig. \ref{fg4} shows BPOs in Dmaps corresponding to BPOs in spatiotemporal patterns in Fig. \ref{fg2}. 
As parameters (CN and Permutation ID) vary, the PBNN exhibits a variety of BPOs. 
The number of CNs (without CN4 and CN6) is $6$ whereas the number of hidden connection parameters $w_{ij}$ is $3^{N^2}$. 
The number of Permutation IDs is $N!$ whereas the number of output connection parameters $c_{ij}$ is $3^{N^2}$. 
In addition, the DBNNs have $2N$ integer threshold parameters $S_i$ and $T_j$.  
It goes without saying that the PBNNs cannot generate more various BPOs than the DBNNs because the PBNNs are included in the DBNNs. 
However, the PBNN parameter space is much smaller than the DBNN parameter space. 
The objective problem is 
{\it relationship between parameters (Permutation ID and CN) and existence/stability of BPOs.}

\section{Globally stable binary periodic orbits}

Depending on parameters, the PBNNs exhibit various BPOs and multiple BPOs can co-exist for initial state. 
Since analysis of multiple BPOs is hard,  
we try to analyze representative BPOs: the globally stable binary periodic orbits (GBPOs). 
This section defines the GBPOs and related concepts. 
First, we note two exceptional endpoints in $\bm{B}^N$:  
\begin{equation}
\bm{x}_- \equiv (-1, \cdots, -1) \in \bm{B}^N, \ \bm{x}_+ \equiv (+1, \cdots, +1) \in \bm{B}^N
\end{equation}
The two endpoints are either fixed points or a BPO with period 2, becuase
\begin{equation}
\begin{array}{ll}
f(\bm{x}_+)=\bm{x}_+, \ f(\bm{x}_-)=\bm{x}_-  \mbox{ if }  w_a+w_b+w_c \ge +1\\
f(\bm{x}_+)=\bm{x}_-, \ f(\bm{x}_-)=\bm{x}_+  \mbox{ if }  w_a+w_b+w_c \le -1
\end{array}
\end{equation}
Hereafter we omit the two endpoints. 
The GBPO is defined by 
\begin{definition}
\label{def2}
A BPO is said to be a globally stable binary periodic orbit (GBPO)
if the BPO is unique (except for $\bm{x}_-$ and $\bm{x}_-$) and 
if all the EPPs fall into the BPO where we assume existence of the EPPs. 
The number of EPPs plus elements of the GBPO is $2^N - 2$. 
\end{definition}
Fig. \ref{fg4}(b) shows a GBPO with period 20 in the Dmap. 
In this 7-dimensional example, $(2^7 - 20 - 2)$ EPPs fall into the GBPO. 
As shown in Section \ref{bfa}, 
depending on the parameters (permutation ID and CN), 
the 7-dimensional PBNNs exhibit a variety of GBPOs and the number of EPPs is more than $2^7/2$. 
The EPPs represent global stability corresponding to error correction \cite{error} of binary signals. 
As the number of EPPs increases, the global stability becomes stronger. 
In the limit case of the M-sequences (e.g., in the linear feedback shift register \cite{lfsr}), 
the period is $2^N$, no EPP exists and is not stable. 
Such M-sequences are different category from the GBPOs in this paper. 
In fundamental viewpoints, uniqueness of the GBPO is convenient to consider existence and stability. 
Analysis of multiple BPOs is complex. 
In application viewpoints, the GBPOs are useful as globally stable signal 
to approximate/predict time-series \cite{uchida} and to control switching circuits \cite{pe1} \cite{pe2} \cite{pe3}.

For simplicity, we focus on the case where $N$ is a prime number $N_p$. 
If an integer $N$ can be factorized into prime factors, classification of the permutation connections becomes complex. 
Here, in order to analyze GBPOs, we define several basic concepts. 

\begin{definition}
\label{def3}
Let $R$ be a shift operator such that 
\begin{equation}
\begin{array}{l}
R: P_0(\sigma_0(1) \cdots \sigma_0(N_p)) \rightarrow P_1(\sigma_1(1) \cdots \sigma_1(N_p))\\
P_1 = R(P_0), \sigma_1(i+1) = \sigma_0(i) + 1 \mbox{ mod }  N_p, i \in \{1, \cdots, N_p\}
\end{array}
\end{equation}
where $\sigma_1(N_p + 1) \equiv \sigma_1(1)$. 
Since the neurons are ring-type connection, 
the permutation connections $P_1$ and $P_0$ ($P$ and $R(P)$) are equivalent even if the permutation IDs are different. 
\end{definition}

\begin{figure}[tb]
  \centering
  \includegraphics[width=0.9\columnwidth]{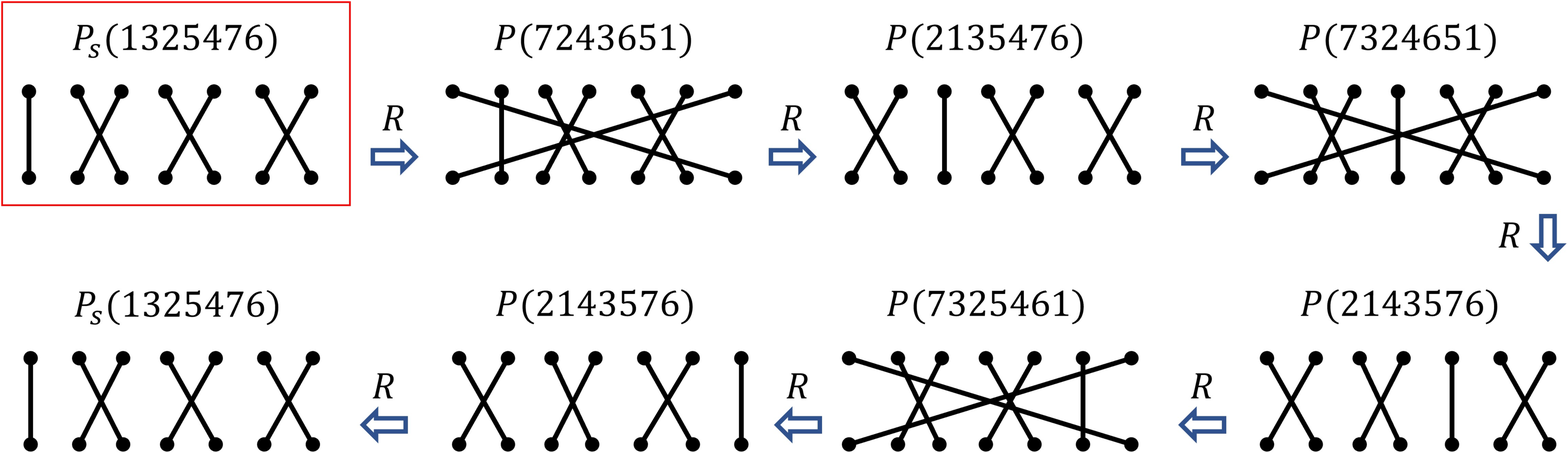}
  \caption{
Equivalent permutation connection examples for $N_p=7$. 
$P_s$: standard permutation connection. $R$: shift operator.
  }
  \label{fg5}

\vspace*{5mm}

  \centering
  \includegraphics[width=0.9\columnwidth]{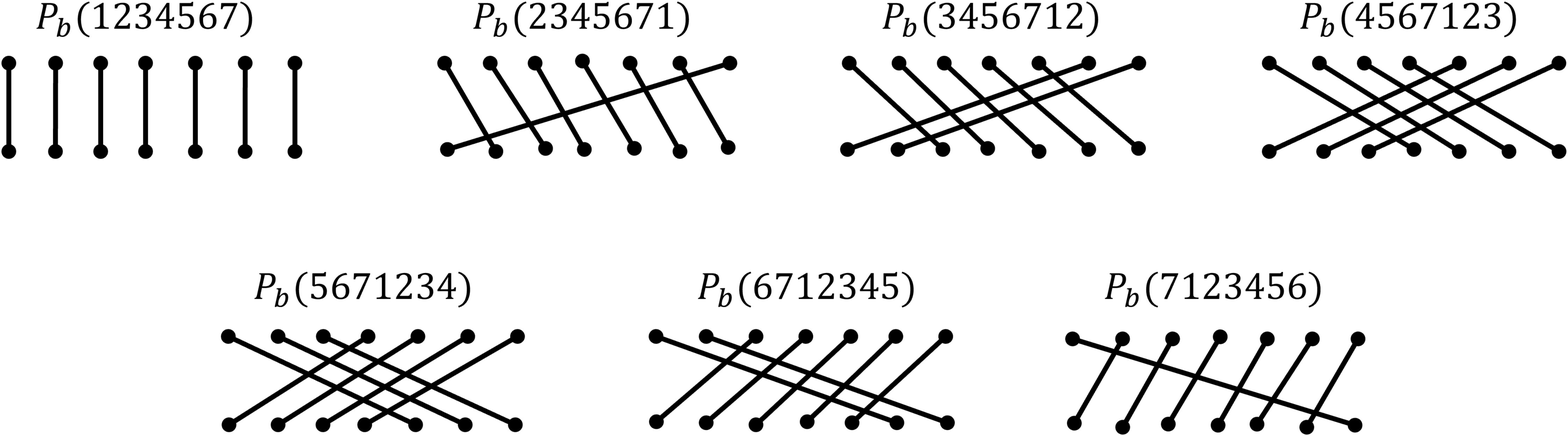}
  \caption{
  7 basic permutation connections for $N_p=7$. 
  }
  \label{fg6}

  \vspace*{5mm}
  
  \centering
  \includegraphics[width=0.66\columnwidth]{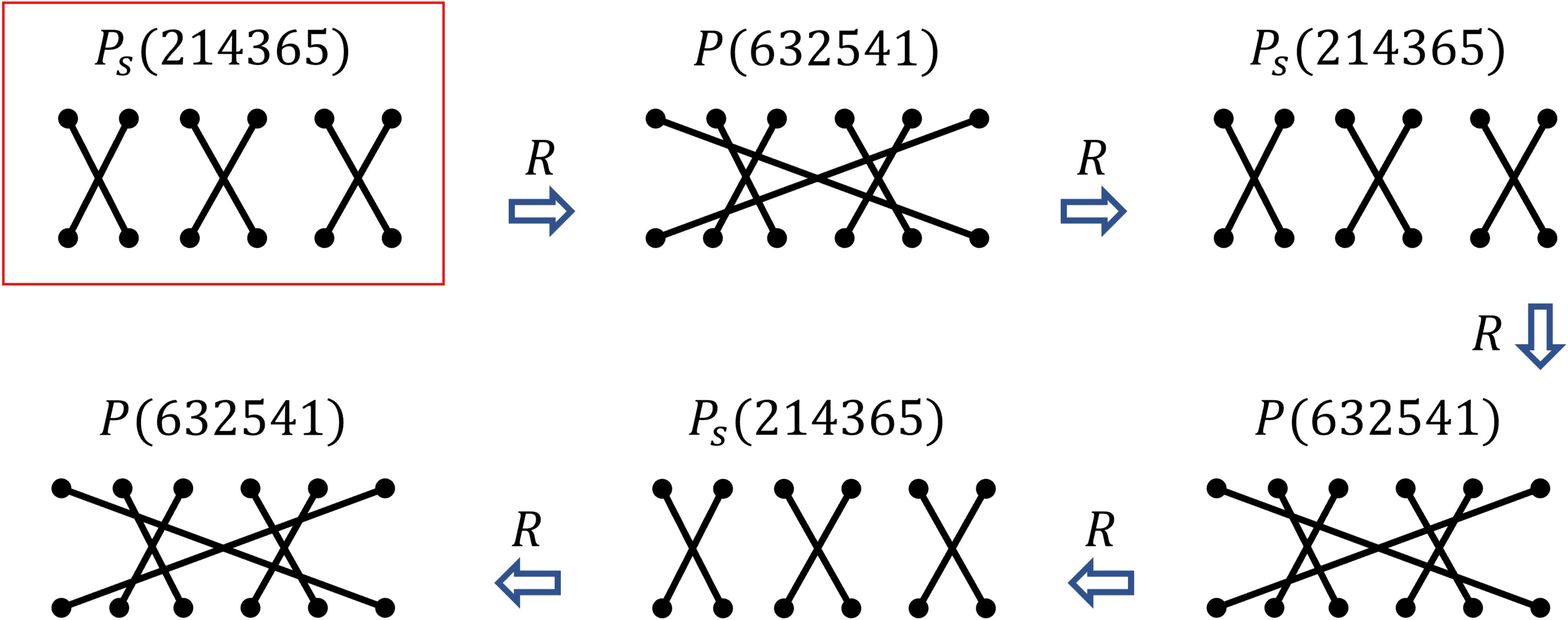}
  \caption{
Permutation connection examples consisting of 3 sub-connections for $N=6$. 
  }
  \label{fg7}
\end{figure}

\begin{definition}
\label{def4}
Let $S$ be a set of permutation IDs that give equivalent permutation connections. 
The set $S$ is referred to as an equivalent permutation set (EPS). 
An EPS is represented by a standard permutation ID $P_s(\sigma_s(0) \cdots \sigma_s(N_p))$
that corresponds to the minimum element in the EPS by means of base-$N_p$ number:
\[
P_s(\sigma_s(1) \cdots  \sigma_s(N_p)) < P_k(\sigma_k(1) \cdots  \sigma_k(N_p)) \in S,  
k \ne s  \mbox{ ( base-$N_p$ number )}
\]
\end{definition}
Fig. \ref{fg5} shows an example of standard permutation connection and 
its equivalent permutation connections for $N_p=7$. 
In this example, the EPS is 
\[
\begin{array}{c}
S= \{ P_s(1325476), P(7243651), P(2135476), P(7324651),\\
 P(2143576),  P(7325461), P(2143576)  \}
\end{array}
\]

\begin{definition}
\label{def5}
A permutation ID $P_b$ is said to be a basic permutation ID  
if it is a fixed point of the shift operator: $R(P_b)=P_b$. 
Since $R(P_b)=P_b$ iff $\sigma_b(i+1) = \sigma{i} + 1$ mod $N_p$, 
the number of basic permutation IDs is $N_p$.  
A basic permutation ID constructs an EPS with one element and is a standard permutation ID. 
\end{definition}

Fig. \ref{fg6} shows basic permutation connections for $N_p=7$. 
Then we have
\begin{theorem}
In $N_p$-dimensional PBNNs, 
the number of standard permutation IDs (i.e., the number of EPSs) is 
$(N_p - 1)! + N_p - 1$ where $N_p \ge 3$ is a prime number. 
\end{theorem}
(Proof) 
Except for $N_p$ basic permutations, 
one standard permutation ID $P_s$ represents $N_p$ equivalent permutation IDs:  
\[
R^{N_p}(P_s) = P_s, \ R^k(P_s) \ne P_s \mbox{ for } 1 \le k \le N_p-1 
\]
where $R^k(P)=R(R^{k-1}(P_s))$ is the $k$-fold composition of the shift operator $R$.  
If there exists an integer $l$ ($2 \le l < N_p$) such that $R^l(P_s) = P_s$, 
the ring-type connection of $P_s$ is decomposed into the same sub-connections  
(e.g., 3 sub-connections $R^{3l}(P_s)=R^{N_p}(P_s) =P_s$ as shown in Fig. \ref{fg7}). 
However, it is impossible for a prime number $N_p$.  
Therefore, except for the basic permutations, the number of standard permutation IDs is $(N_p! - N_p)/N_p$. 
Adding the $N_p$ basic permutation IDs, the number of standard permutation IDs is 
$(N! - N_p)/N_p + N_p = (N_p - 1)! + N_p -1$.

\clearpage

\section{Brute force attack to explore GBPOs}
\label{bfa}
Table \ref{tb0} shows the number of standard permutation IDs for prime numbers $N_p$ 
together with the number of full binary connection parameters between hidden and output layers in the DBNNs for $N=M=N_p$. 
The number of the permutation connections is much smaller than the number of the full binary connections. 
However, analysis of the GBPOs becomes harder as $N_p$ increases. 
For convenience, we consider GBPOs in 7-dimensional PBNNs ($N_p=7$). 
In the case $N_P=7$, 
the number of all the standard permutation connections is $(N_p-1)! + N_p-1 = 726$, 
the number of initial points is $2^7$, and the brute force attack is possible.  
We can clarify the number and period of all the GBPOs precisely. 
Analysis of the 7-dimensional GBPOs are fundamental to consider higher-dimensional GBPOs 
and their engineering applications. 

\begin{table}[b!]
  \centering
  \caption{The number of standard permutation connections in PBNN 
  and full binary connections between hidden and output layers in DBNN.}
  \label{tb0}
  \begin{tabular}{|c|c|c|} \hline
  $N_p$  &  \# standard permutation IDs  & \# full binary connections \\ \hline
  3   &  4                & $2^{9}$  \\
  5   &  28               & $2^{25}$ \\
  7   &  726              & $2^{49}$ \\
  11  &  3628810          & $2^{121}$\\
  13  &  479001612        & $2^{169}$\\
  17  &  20922789888016   & $2^{289}$\\ \hline
  \end{tabular}
  \end{table}  

We explore the 7-dimensional GBPOs as the following. 
First, as state earlier, objective connection numbers are 
CN0, CN2, CN2, CN3, CN5, and CN7 (CN1 $\equiv$ CN4 and CN3 $\equiv$ CN6).   
Second, applying the shift operator $R$, we obtain the 726 standard permutation IDs. 
Third, applying the brute force attack to each standard permutation ID and CN, 
we obtain BPOs and their EPPs where we use the BPO calculation algorithm in \cite{horimoto}. 
If the number of a BPP plus its EPPs is $2^7-2=126$ then the BPO is declared as the GBPO. 
The period of the GBPO is stored together with its standard permutation ID. 

\begin{figure}[htb]
  \centering
  \includegraphics[width=0.8\columnwidth]{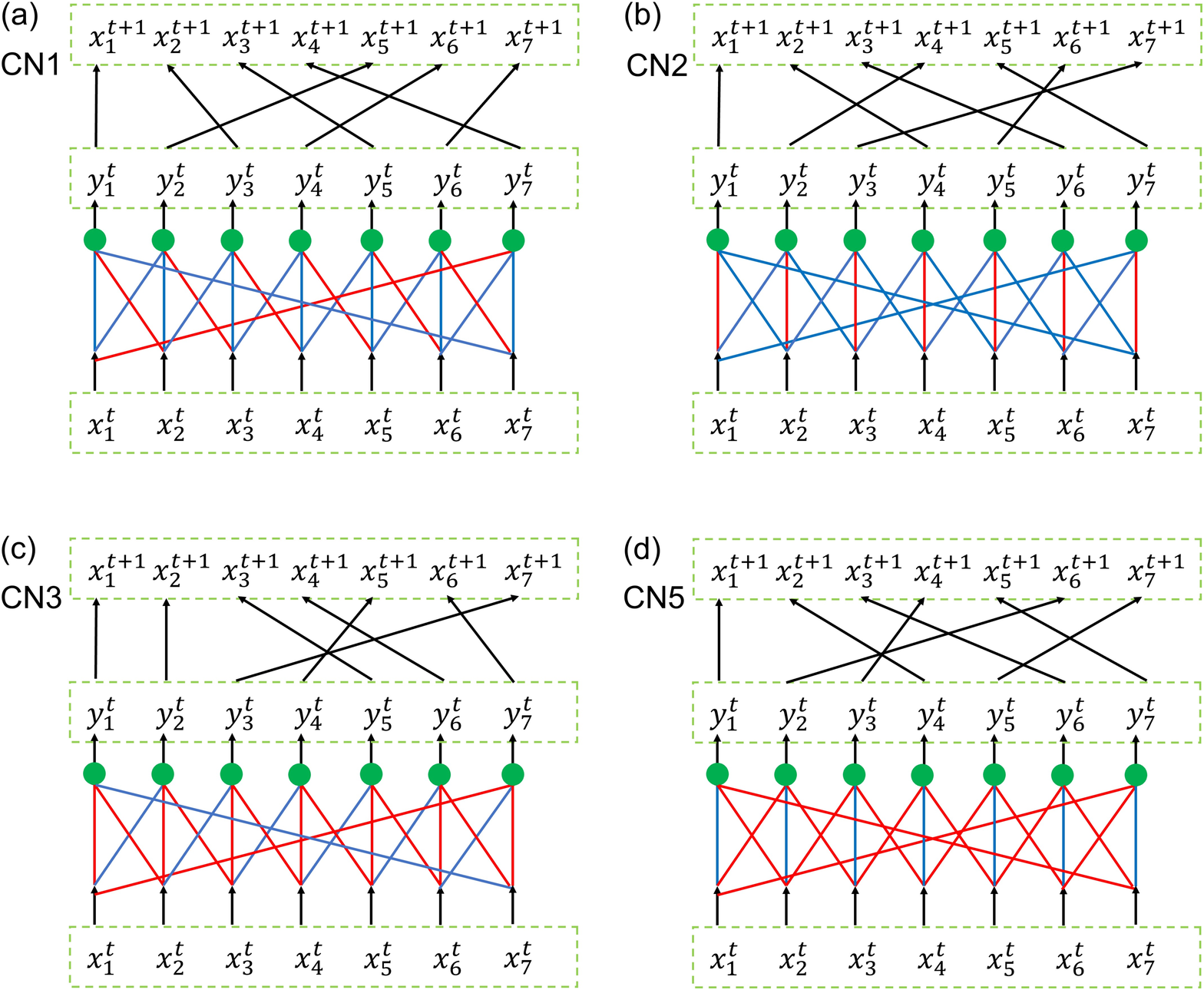}
  \caption{PBNN examples (exhibit GBPOs) for $N_p=7$. 
(a) $P_s(1357246)$, CN1. 
(b) $P_s(1462753)$, CN2. 
(c) $P_s(1256473)$, CN3.
(d) $P_s(1463725)$, CN5. 
  }
  \label{fg8}
  \end{figure}

\begin{figure}[htb]
\centering
\includegraphics[width=1.0\columnwidth]{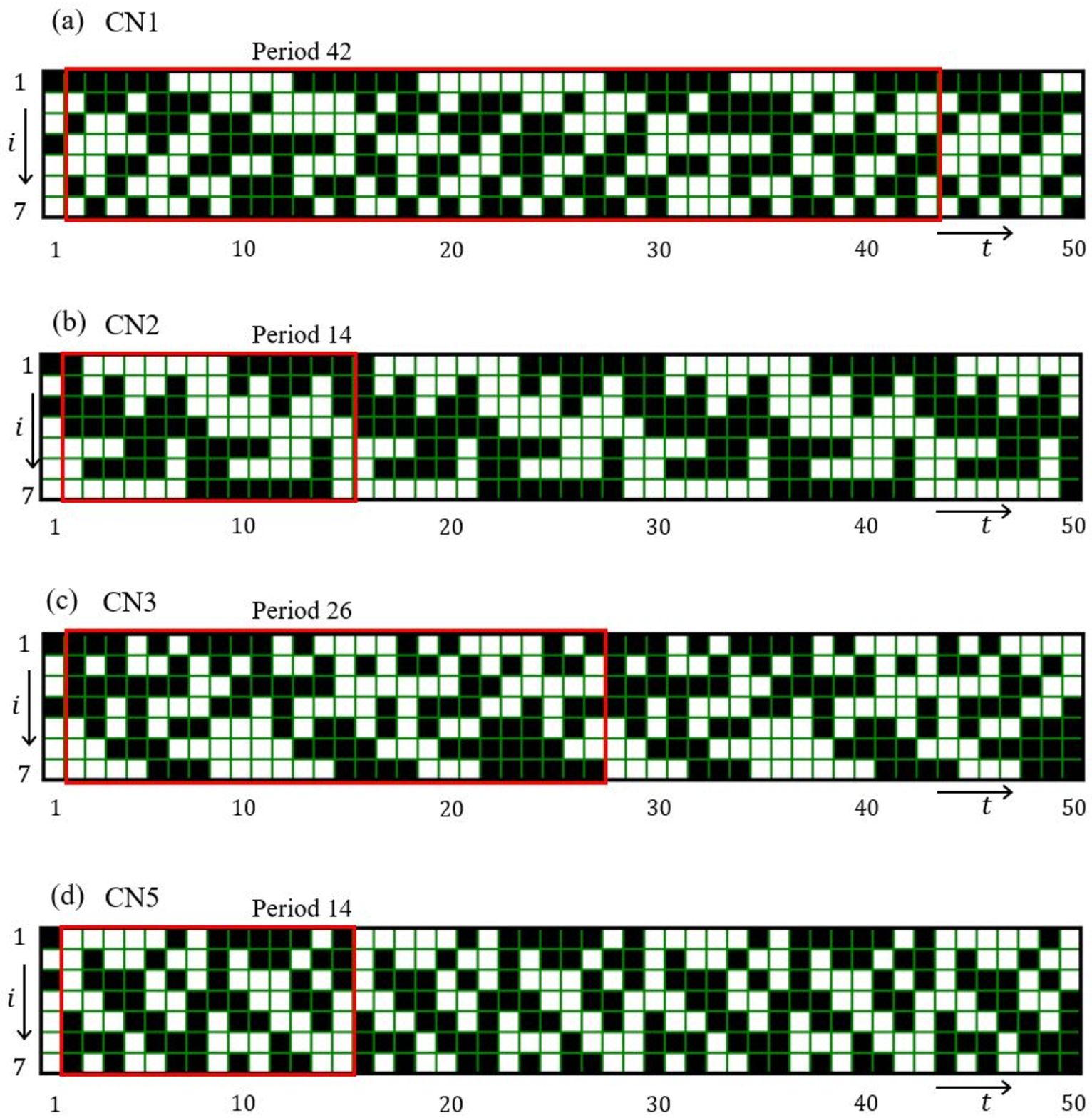}
\caption{GBPO examples as spatiotemporal patterns for $N_p=7$. 
(a) $P_s(1357246)$, CN1, GBPO with period 42. 
(b) $P_s(1462753)$, CN2, GBPO with period 14. 
(c) $P_s(1256473)$, CN3, GBPO with period 26.
(d) $P_s(1463725)$, CN5, GBPO with period 14. 
}
\label{fg9}
\end{figure}

\begin{figure}[htb]
\centering
\includegraphics[width=1.0\columnwidth]{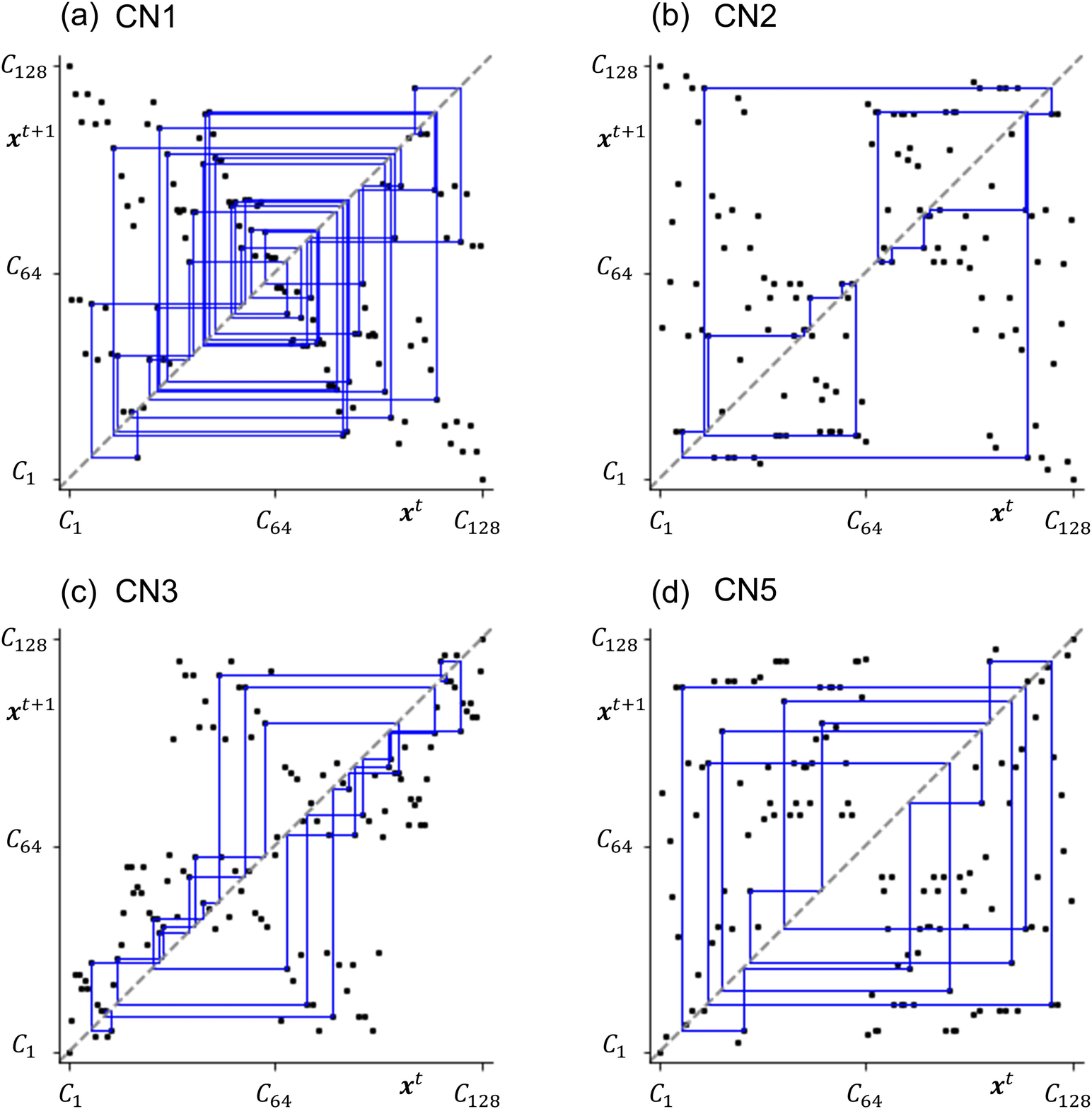}
\caption{GBPO examples in Dmaps. 
(a) $P_s(1357246)$, CN1, GBPO with period 42. 
(b) $P_s(1462753)$, CN2, GBPO with period 14. 
(c) $P_s(1256473)$, CN3, GBPO with period 26.
(d) $P_s(1463725)$, CN5, GBPO with period 14. 
}
\label{fg10}
\end{figure}

In the exploration, it is confirmed that CN0 and CN7 cannot provide GBPO.   
The CN0 and CN7 are omitted hereafter. 
Fig. \ref{fg8} shows typical examples of PBNNs for CN1, CN2, CN3, and CN5 
that generate GBPO with period 42, period 14, period 26, and period 14, respectively. 
Fig. \ref{fg9} shows the 4 GBPOs as spatiotemporal patterns and 
Fig. \ref{fg10} shows the 4 GBPOs in Dmaps. 

As a criterion of the period, we give 
\begin{definition}
\label{def6}
For identity permutation ($P_b(1234567)$ for $N_P=7$), 
the period of the BPO is said to be basic period. 
If the PBNN generates multiple BPOs, the maximum period is adopted.  
\end{definition}
For CN1, the basic period is 14 as a BPO in Fig. \ref{fg2} (a) that is a GBPO. 
We have confirmed that the identity permutation $P_b(1234567)$ cannot provide a GBPO.  
In Figs. \ref{fg8} to \ref{fg10}, we can see that, adjusting permutation IDs from the identity permutation $P_b(1234567)$, 
the PBNNs can generate a variety of BPOs represented by the GBPOs with longer period. 
As the main result, 
tables \ref{tb1} to \ref{tb4} show 
a list of standard permutation IDs for all the GBPOs.  
As stated in Definition \ref{def4}, 
each standard permutation ID represents 7 equivalent permutation IDs. 
We give an overview of the list for CN1, CN2, CN3, and CN5:  
\begin{itemize} 
\item CN1: The basic period is 14 for $P_b(1234567)$. 
The PBNNs generate 27 GBPOs. 
The maximum period is 42 for $P_s(1357246)$ as shown in Fig. \ref{fg9} (a). 
The number of EPPs is $126-42$. 

\item CN2: The basic period is 2. 
The PBNNs generate 56 GBPOs. 
The maximum period is 14 where the number of EPPs is $126-14$, 
e.g. $P_s((1462753)$ as shown in Fig. \ref{fg9} (b). 

\item CN3: The basic period is 14. 
The PBNNs generate 28 GBPOs. 
The maximum period is 26 where the number of EPPs is $126-26$, 
e.g. $P_s(1256473)$ as shown in Fig. \ref{fg9} (c). 

\item CN5: The basic period is 2. 
The PBNNs generate 62 GBPOs. 
The maximum period is 14 where the number of EPPs is $126-14$, 
e.g. $P_s(1463725)$ as shown in Fig. \ref{fg9} (d).
 
\end{itemize} 
These tables clarify relation between parameters (permutation ID and CN) and periods of the GBPOs. 
The number of EPPs is 126 minus the period. 
As the parameters vary, the 7-dimensional PBNNs can generate a variety of GBPOs.  
These results provide fundamental information to analyze various PBNNs and to synthesize PBNNs with desired GBPOs.

\clearpage

\begin{table}[thb]
\centering
\caption{Standard permutation ID and period of GBPO for CN1}
\label{tb1}
\begin{tabular}{|cc|cc|cc|} \hline
    ID   & period   &   ID  &  period &   ID  &  period \\ \hline
1256374  &  26   & 1625473 &  6 & 2517436 & 18  \\
1257436  &  18   & 1627435 & 16 & 2576314 & 12  \\
1273654  &  14   & 1657234 & 12 & 2613754 & 20  \\
1352476  &  34   & 1657243 &  4 & 2615374 & 12  \\
\bf{1357246}  &  \bf{42} & 1672453 & 18 & 2675314 & 8  \\
1375426  &  26   & 1672543 &  2 & 2751436 & 8  \\
1526374  &  42   & 1673425 & 18 & 2763154 & 20  \\
1527643  &  14   & 2175346 & 10 & 3416725 &  8  \\
1576324  &  24   & 2417356 & 14 & 4671325 & 24  \\ \hline
\end{tabular}
\end{table}

\begin{table}[thb]
  \centering
  \caption{Standard permutation ID and period of GBPO for CN2}
  \label{tb2}
  \begin{tabular}{|cc|cc|cc|} \hline
   ID  & period   &   ID  &  period &   ID  &  period \\ \hline
  1367245  &  4   & 1653724 & 14 & 2641735 & 2   \\
  1427365  &  8   & 1674352 & 6  & 2641753 & 2   \\
  1436275  &  8   & 1675234 & 4  & 2671354 & 2   \\
  1436752  &  8   & 1732645 & 8  & 2671453 & 2   \\
  1457236  &  4   & 1742653 & 2  & 2761345 & 4   \\
 \bf{1462753}  &  \bf{14}  & 1762453 & 2  & 2761354 & 2   \\
  1467325  &  6   & 1762543 & 6  & 2761453 & 2   \\
  1476235  &  6   & 2156734 & 14 & 2763154 & 2  \\
  1476532  &  6   & 2167345 & 14 & 3156724 & 4   \\
  1527643  &  10  & 2356714 & 8  & 3176254 & 4   \\
  1564372  &  6   & 2361754 &  2 & 3471256 & 4   \\
  1567423  &  4   & 2365174 & 2  & 3517426 & 2   \\
  1572643  &  14  & 2367145 & 4  & 3571246 & 2  \\
  1627543  &  2   & 2367154 &  2 & 3576214 & 4   \\
  1642735  &  14  & 2461753 & 2  & 4167253 & 2   \\
  1643752  &  2   & 2467135 &  2 & 4173256 & 2   \\
  1645732  &  2   & 2467153 & 2  & 4617325 & 4   \\
  1647532  &  2   & 2517643 & 4  & 4712356 & 8  \\
  1653274  &  10  & 2571634 & 2  &         &     \\ \hline
  \end{tabular}
  \end{table}

  \begin{table}[thb]
    \centering
    \caption{Standard permutation ID and period of GBPO for CN3}
    \label{tb3}
    \begin{tabular}{|cc|cc|cc|} \hline
    ID  & period   &   ID  &  period &  ID  &  period \\ \hline
    1235476  &  14   & 1567243 & 12 & 2761345 & 2   \\
    1246753  &  22   & 1576324 & 24 & 3157426 & 12  \\
    \bf{1256473}  &  \bf{26}   & 1652473 & 10 & 3167425 & 8   \\
    1267435  &  26   & 1657243 &  4 & 3176245 & 10  \\
    1362754  &  6    & 2156374 & 22 & 3561724 & 10  \\
    1375462  &  2    & 2417635 &  6 & 3567214 & 24  \\
    1425376  &  10   & 2463175 & 20 & 3612745 & 8   \\
    1463275  &  6    & 2516374 & 2  & 3761425 & 12  \\
    1465273  &  16   & 2516473 & 8  &         &     \\
    1476235  &  2    & 2641753 & 20 &         &     \\ \hline
    \end{tabular}
    \end{table}
 
    \begin{table}[thb]
      \centering
      \caption{Standard permutation ID and period of GBPO for CN5}
      \label{tb4}
      \begin{tabular}{|cc|cc|cc|} \hline
      ID  & period   &   ID  &  period &   ID  &  period \\ \hline
      1245673 &	10 & 1436752	& 8	 & 1673452	&10 \\
      1246735	& 2	& 1437526	& 2	 	& 1675234	&4  \\
      1247635	& 6	& 1453726	& 6	 	& 1726543	&2  \\
      1256734	& 10 	& 1457236	& 4	 	& 1732645	&8  \\
      1257346	& 2	& 1457263	& 10 		& 1745326	&6  \\
      1257436	& 6	& 1463275	& 6		& 1745623	&10 \\
      1267345	& 10 	& \bf{1463725}	& \bf{14} 	& 1756243	&2  \\
      1273456	& 10 	& 1465723	& 10 		& 1765243	&6  \\
      1342765	& 6  	& 1467235	& 2	 	& 2356714	&8  \\
      1356724	& 10 	& 1467523	& 10 		& 2367145	&4  \\
      1367245	& 4	& 1472635	& 14 		& 2517643	&4  \\
      1372645	& 6	& 1532764	& 2		& 2761345	&4  \\
      1374265	& 2	& 1547326	& 6		& 3156724	&4  \\
      1376254	& 6	& 1567423	& 4	 	& 3176254	&4  \\
      1376452	& 6	& 1572346	& 10 		& 3461725	&10 \\
      1423765	& 6	& 1572364	& 10 		& 3471256	&4  \\
      1427365	& 8	& 1654372	& 2	 	& 3576214	&4  \\
      1427635	& 6	& 1657342	& 2	 	& 3617245	&10 \\
      1432675	& 6	& 1657432	& 6	 	& 4617325	&4  \\
      1432756	& 6	& 1672435	& 10 		& 4712356	&8  \\
      1436275	& 8	& 1672534	& 10 		& 	        & \\ \hline
      \end{tabular}
      \end{table}

\clearpage

\section{Conclusions}
Fundamental dynamics of the PBNNs has been studied in this paper. 
The PBNNs are characterized by global permutation connections and local binary connections. 
Although the parameter space is much smaller than existing recurrent-type neural networks, 
the PBNN can exhibit various BPOs.  
In order to realize precise analysis, we focus on the GBPOs and define standard permutation connections. 
Applying the brute force attack to 7-dimensional PBNNs, 
we have presented complete list that clarifies relationship between parameters and periods of GBPOs. 
Even in the 7-dimensional cases, the PBNNs exhibit a variety of GBPOs. 
It suggests that higher dimensional PBNNs exhibit a huge variety of BPOs/EPPs. 
Many problems remain in our future works: 

\begin{itemize}

\item Mechanism to generate the GBPOs. 

\item Classification and stability analysis of various BPOs. 
Besides the GBPOs, the PBNNs exhibit various BPOs, depending on parameters and initial conditions. 
 
\item Effective evolutionary algorithms \cite{mop1} \cite{mop2} for 
analysis of higher dimensional BPOs where the brute force attack is impossible. 

\item Effective evolutionary algorithms for 
synthesis of PBNNs with desired BPOs. 

\item Efficient hardware implementation for engineering applications including 
robust control signals of switching circuits and time-series approximation/prediction. 
The PBNNs are well suited for FPGA based hardware implementation 
that transforms the BPOs into electric signals in the applications.

\end{itemize}

\section*{Declaration of competing interest}
The authors declares that he has no known competing financial interests or personal
relationships that could have appeared to influence the work reported in this paper.


\end{document}